\documentclass[11pt]{amsart}
\usepackage[latin1]{inputenc}
\usepackage[cyr]{aeguill}
\usepackage{amsxtra}
\usepackage{amssymb}
\usepackage[all]{xy}
 
\newtheorem{theorem}{\sc Theorem}
\newtheorem{proposition}{\sc Proposition}
\newtheorem*{lemma}{\sc Lemma}
\newtheorem{corollary}{\sc Corollary}
\theoremstyle{remark}
\newtheorem{remark}{\it Remark}
\newtheorem{example}{\it Example}
\newtheorem{definition}{\sc Definition}

\font\tmsb=msbm10 at12pt
\font\smsb=msbm7
\font\ssmsb=msbm5
\newfam\msbfam
\textfont\msbfam=\tmsb
\scriptfont\msbfam=\smsb
\scriptscriptfont\msbfam=\ssmsb
\def \RM{\mathbb {R}}
\def \ZM{\mathbb{Z}}
\def \CM{\mathbb{C}}

\def \Gt {\mathcal{G}}

\def \d{\partial}

\def\dt{\delta} 
\def\a{\alpha}
\def\b{\beta}
\def\e{\varepsilon}  
\def\g{\gamma}
\def\p{\varphi}

\def\l{\lambda}

\def\lb{\left\{}
\def\rb{\right\}}
\def\G{\Gamma}   
\def \S{\Sigma}
\def \t{\tilde}

\def \to{\longrightarrow} 
\def \w{\wedge}
\def\2n{(\CM^{2n},0)}

\def \gk{\mathfrak{g}}
\def \hk{\mathfrak{h}}

\def \< {{\langle }}
\def \> {{\rangle }}
\def \Ker {{\rm Ker}}
\def \Hom {{\rm Hom}}
\def \sl {{\mathfrak sl}}

\newcommand{\OM}{{\mathcal O}}

\newcommand{\D}{\Delta }

\medskip
\begin{document}
\title [Vanishing cycles in complex symplectic geometry]
{Vanishing cycles in \\
complex symplectic geometry}
\author[Mauricio D. Garay]{ Mauricio D. Garay$^*$  }
\date{December 2005}
\thanks{$*$ partly supported by the Forschungsstipendium GA 786/1-1
of the Deutsche Forschungsgemeinschaft and by the
IH{\'E}S ($6^{th}$ European Framework program,
contract Nr. RITA-CT-2004-505493)}
\date{October 2006}
\address{SISSA/ISAS, via Beirut 4, 34014 Trieste,
Italy.}
\email{garay@sissa.it}
\thanks{\footnotesize 2000 {\it Mathematics Subject Classification:}
32S50}
\keywords{Monodromy, vanishing cycles, integrable systems,
Symplectic geometry, Lagrangian varieties,
involutive varieties, Simple Lie algebras.}
\begin{abstract}
We study the vanishing cycles on the Milnor fibre for some non-isolated singularities which appear naturally
in symplectic geometry. Under
assumptions given in the text, we show that
the vanishing cycles associated to a distinguish basis
freely generate the corresponding homology groups of the Milnor fibre.
We derive some consequences of this fact, in particular for the study of adjoint orbits in Lie algebras.
\end{abstract}
\maketitle
\parindent=0cm
\section*{Introduction}
The local Picard-Lefschetz theory studies the vanishing cycles
associated to a holomorphic function germ
$ f:(\CM^m,0) \to (\CM,0)$ (\cite{Lefschetz,Mil,Pham_monodromie}). The first case to be studied has been that of isolated
singularities, then using the Whitney-Thom theory of stratified spaces and morphisms,
Deligne, Grothendieck, Hironaka, L\^e and others were able to develop a general theory 
which includes the case of non-isolated singularities
(\cite{Deligne_cycles,Sabbah_these,Hironaka_oslo,Le_fourier}).\\
Using part of these works, we show that the non-isolated
singularities which appear in symplectic geometry, as for instance
in finite dimensional integrable systems, behave to a large extent like
isolated hypersurface singularities.
To this aim, we shall make an axiomatisation of these type of non-isolated singularities.\\
Our main result states that, under assumptions given in the text,
the Lefschetz type vanishing cycles freely generate the corresponding homology group.
This fact, which is already wrong for arbitrary isolated complete intersection singularities, shows that
the non-isolated complete intersection singularities that we consider are, from the point of view of topology, the simplest one
after the case of isolated hypersurface singularities.\\
Indeed, our approach is purely topological and elementary, close in spirit to the original works of
Brieskorn, L\^e, Milnor and Pham; it does not use the Deligne-Grothendieck
formalism nor any sheaf theoretic methods. The philosophy is in some sense opposite: by elementary topological considerations,
we hope deduce non-trivial facts on deformation theory and quantisation; this is however at the present stage only speculative.
Although the applications outside the realm of differential topology are not firmly established, we can apply
our theorem for investigating the topology of non-isolated singularities;
we shall give two examples: involutive mappings in symplectic spaces (which includes integrable systems in symplectic spaces) and the Steinberg mapping associated to a semi-simple Lie algebra. The case of integrable systems
is particularly relevant for perturbation theory since we
now that the Lagrangian deformation theory is essentially
topological (\cite{mutau}).
 Most of
the results that we obtained are part of the following program
formulated by Arnold at his seminar (\cite{Arnold_Problem}, 1992-2): ``(...) Construct a monodromy and variation
theory for non-isolated singularities (...). The obtained theory must apply to the mapping
that assigns characteristic polynomials to matrices. It must generalize the
Brieskorn-Grothendieck description of simple singularities for the A,D,E surfaces over non
quasi-regular elements of Lie algebras.''
\section{Non-isolated complete intersection singularities}
We start by showing how non-isolated singularities typically appear in symplectic geometry; then we review classical facts
on non-isolated complete intersection singularities.
\subsection{Propagation of singularities in symplectic geometry}
Let us consider the space $\CM^{2n}=\{ (q,p) \}$ together with the symplectic form
$\sum_{i=1}^n dq_i \w dp_i$ and a holomorphic mapping $f=(f_1,\dots,f_k):M \to \CM^k$
whose components  are in involution ,i.e.,
$$\forall i,j \leq k,\ \  \lb f_i,f_j \rb =\sum_{l=1}^n \d_{q_l} f_i \d_{p_l}f_j- \d_{p_l} f_i
  \d_{q_l}f_j=0. $$
Here $M \subset \CM^{2n}$, $S \subset \CM^n$ denote open subsets.\\ 
The fibres of such an involutive mapping might be, for instance, characteristic varieties of
partial differential equations depending on a semi-classical
parameter.\\
The Hamiltonian vector fields associated to the components of $f$ are tangent to its fibres,
therefore such an involutive map may generically have fibres
with a non-isolated singular locus.\\
As a simple example, consider the map $f=(f_1,f_2):\CM^4 \to \CM^2$ defined by the equations
$f_1=p_1q_1$, $f_2=p_2$ and denote by $V_0$ its fibre above the origin.
The variety $V_0$ is analytically equivalent to the product of the plane curve
$\G=V_0 \cap \{ q_2=0 \}$ by an affine line.
The singular locus of the variety $V_0$ is one
 dimensional and the Hamiltonian vector field $X=\d_{q_2}$ associated to $f_2$ is tangent to it.\\
This example is stable in the sense
 that deforming the germ of $f$ at the origin by a map whose components are in
 involution, leaves the analytical type of the singular fibration germ
 given by $f$ unchanged (\cite{lagrange,moment}).\\
Remark that the smooth fibres of $f$ have a vanishing topology which
is drastically different from that of isolated complete
intersection singularities.
In such a case, by the classical results of Milnor and Hamm (\cite{Hamm,Mil}),
the Milnor fibre of the germ $f:(\CM^4,0) \to (\CM^2,0)$ would have been
homeomorphic to a bouquet of 2-dimensional spheres while it is homeomorphic
to a circle.\\
This is a general feature for singular Lagrangian varieties since the first Betti number of the fibre is equal to the dimension of the Lagrangian versal deformation space (\cite{mutau}).

\subsection{The Milnor fibration for non-isolated complete intersection singularities}
Consider a morphism $f:X \to S$ between complex spaces over a smooth base.
There is, a priori, no reason for $f$ to define a topologically trivial fibration above its regular values,
even if $X$ is smooth; Thom's $a_f$-condition ensures such a result.
\begin{definition}A continuous map between Whitney stratified topological spaces $f:X \to S$, $X=\cup_{\a}X_\a$
, $S=\cup_{\a} S_\a$
is called {\em stratified} if
the restriction of $f$ to each strata is a submersion.
\end{definition}
For a given holomorphic map $f:X \to S$, the union of singular loci of its fibres 
form an analytic space called the {\em critical locus} of the map, its image under the map is called the {\em discriminant}.
The discriminant is a priori only a variety and does not carry a canonical scheme structure
(it is nevertheless very likely that for the particular cases we will study, it does carry such a structure analogous to that discovered in \cite{Teissier_fitting} for isolated singularities).
\begin{definition}[\cite{Thom_strates}]
A stratified map $f:X \to S,\ X \subset \RM^n$ satisfies the {\em condition $a_f$}
if for any sequence of point $(x_n)$ in a strata $X_{\a}$ converging to a point $x$ lying in an
adjacent strata $X_{\a'}$ such that the sequence $\Ker\, df_{|X_{\a}}(x_n)$ converges, we have the inclusion
$\Ker\, df_{|X_{\a'}}(x) \subset \underset{n \to \infty}{\lim} \Ker\, df_{|X_{\a}}(x_n)$.
 The map $f$ satisfies the {\em strong $a_f$-condition} if in addition the complement of the critical
set of $f$ is a strata.
\end{definition}
A function germ satisfies the $a_f$-condition if it admits a representative satisfying this condition.
Denote by $B_\dt \subset \CM^m$, the closed ball centred at the origin of radius $\delta$.
\begin{definition}A {\em standard representative $\bar f:X \to S$} of a holomorphic map germ $f:(Y,0) \to (\CM^k,0)$, $Y \subset \CM^m$,
is a proper representative which satisfies the following conditions
\begin{enumerate}
\item there exists a surjective holomorphic map $g:Y \to S'$ and a compact neighbourhood $S \subset S'$
of the origin so that $\bar f$ is the restriction of $g$ to $X= g^{-1}(S) \cap B_{\delta}$,
for some closed ball $B_\delta$ contained in $Y$,
\item the map $g$ satisfies the $a_f$-condition, 
\item the fibres of $g$ are transversal to the boundary of the ball $B_\delta$,
\item the spheres centred at the origin of radius $r \leq \delta$ are transversal to the fibre of $g$ above the origin.
\end{enumerate}
\end{definition}
In the sequel, we shall abusively use the same letter for the germ and a standard representative of it.
A holomorphic map germ $f:(\CM^m,0) \to
(\CM^k,0)$ whose special fibre is the germ of a reduced complete intersection will be called a
{\em complete intersection map germ}.\\
A small modification of Ehresmann's lemma gives the following result.
\begin{proposition}
Consider a standard representative $f:X \to S$ of a complete intersection map germ.satisfying the $a_f$-condition.
The restriction of $f$ above the open stratum defines a $C^\infty$ locally trivial fibration.
\end{proposition}
\begin{example}
Denote by $B$ the closed unit ball in $\CM^4$.
Put $X=\{ (x,y,z,\l) \in B:(xy+\l)z=1$ and consider the projection on the last factor
$$f:X \to S,\  (x,y,z,\l) \mapsto \l.$$
We consider the Whitney stratification with only two strata: one
is the interior of $X$, the other one is the boundary of $X$.
The fibre of $f$ above $\l$ is diffeomorphic to the complement of the curve $\{ xy=\l \}$ in the affine plane
$\CM^2=\{ (x,y) \}$; therefore for $\l=0$ it is homotopic to a 2-torus. It is a simple exercise to see that the
first Betti number of the other fibres equals 1  (one can also prove that they are actually homotopic to the complement of a circle in $\RM^3$).
Therefore the stratified map $f:X \to S$ does not satisfy the strong $a_f$-condition.
\end{example}
The fibration associated by a standard representative of a holomorphic
map germ $f:(\CM^{2n},0) \to (\CM^k,0)$ does not depend on the choice of the representative, it is called the {\em Milnor fibration} of $f$. The fibre of a standard representative above the origin is called the {\em special fibre}.
\subsection{The local Lefschetz theorem}
\label{SS::Lefschetz}
\begin{theorem}[\cite{Le_fourier}]
\label{T::Le}
Let $f:(X,0) \to (\CM,0), X \subset \CM^m,$
be the germ of a complex analytic map satisfying the $a_f$-condition. For any smooth
Milnor fibre $V$, there exists
an Zariski open subset $\Omega$ in the space of linear hyperplanes $H \subset \CM^m$ such that
the pair $(V,H \cap V)$ is $N$-connected with $N=dim V-1$.
 \end{theorem}
\begin{remark}The theorem was originally stated for a smooth manifold $X$ but this condition is used only
to ensure that a standard representative of $f$ satisfies the $a_f$-condition.
\end{remark}

\begin{theorem}[\cite{Le_oslo,Hironaka_oslo}]
\label{T::Hironaka}
A holomorphic function germ $f:(X,0) \to (\CM,0)$ defined on a reduced complex space $X$
satisfies the $a_f$-condition.
\end{theorem}

Using the two above-mentioned theorems, we get the following generalisation of a theorem due to Kato-Matsumoto
\cite{Kato_Matsumoto}.
\begin{corollary}
\label{C::Lefschetz}
Consider a complete intersection map germ $f:(\CM^m,0) \to
(\CM^k,0)$ satisfying the $a_f$-condition and denote by $s$ the dimension of the singular locus
of the special fibre.
The smooth fibres of a standard representative of $f:X \to S$ are $(m-k-s-1)$-connected. 
\end{corollary}
\begin{proof}
The proof is completely similar to the one given by Ch\'eniot and L\^e for hypersurface singularities \cite{Le_Cheniot}.\\
We use induction on the dimension $s$ of the singular locus of the special fibre.
For $s=0$, the result is due to Milnor and Hamm (\cite{Hamm,Mil}).\\ Assume that $s>0$,
the subset $\Omega \subset \Hom_{\CM}(\CM^k,\CM^{k-1})$ for which the zero
fibre $X_{\pi}$ of the map $\pi \circ f:(\CM^m,0) \to (\CM^{k-1},0)$ is a reduced complete intersection germ,  is readily seen to be a Zariski open subset (in case $k=1$, we take $X_{\pi}=\CM^m$).\\
Take $\pi \in \Omega$ and chose $f_j$ so that the ideal generated by the components of $f$ is the same than that generated by the components of $\pi \circ f$ and $f_j$.
Denote by $g:X_{\pi} \to (\CM,0)$ the restriction of the function $f_j$ to $X_{\pi}$. Theorem \ref{T::Hironaka}
implies that for a well chosen stratification, the map $g$ satisfies the $a_f$-condition.\\
Theorem \ref{T::Le} implies in turn that
for a generic hyperplane $H$ and any smooth fibre $V$ associated to a standard representative of
$g$, the pair $(V,V \cap H)$ is $(m-k-1)$-connected. The induction assumption states that
$V \cap H$ is $(m-k-s-1)$-connected, therefore $V$ is also  $(m-k-s-1)$-connected. This concludes the proof of the corollary.
\end{proof}

\section{The Milnor fibration for $T_n$-type maps}
\subsection{Transverse maps}
Recall that
\begin{enumerate}
\item two continuous map germs $f,g:(\CM^m,0) \to (\CM^k,0)$ are called
{\em topologically equivalent} if there exists a commutative diagram
$$
\xymatrix{
(\CM^m,0) \ar[r]^f \ar[d] &  (\CM^k,0) \ar[d]\\
(\CM^m,0) \ar[r]^g & (\CM^k,0)
}
$$
where the vertical arrows are homeomorphisms.
\item a holomorphic map-germ $f:(\CM^n,0) \to (\CM^k,0)$ belongs to the {\em $A_1$ singularity class}
if its quadratic differential at the origin is definite.
\end{enumerate}
The complex Morse lemma states that a map which belongs to the $A_1$ singularity class
can always be reduced to a sum of squares after a holomorphic change of coordinates.
\begin{definition}
A {\em suspension} of a map germ $f:(\CM^n,0) \to (\CM^k,0)$ is a map germ of the type
$$\t f:(\CM^n \times \CM^l,0) \to (\CM^k \times \CM^s,0),\ (x,y) \mapsto (f(x),\pi(y))$$
where $\pi$ is a linear surjective mapping.
\end{definition} 
The critical values for which the germ of $f$ is topologically equivalent to
a suspension of an $A_1$ singularity $g:(\CM^{n+1},0) \to (\CM,0)$ in $(n+1)$ variables
will be called {\em $ TA_1^n $ critical values}.
\begin{definition}
A complete intersection map germ $f:(\CM^m,0) \to (\CM^k,0)$ is called
of {\em $T_n$-type} if it satisfies the strong $a_f$-condition and
if for any standard representative $f:X \to S$,
there is an open dense  subset of $TA_1^n$ critical values in the discriminant for which
the singular locus of the variety $f^{-1}(z)$ is connected.
\end{definition} 
\begin{example}
Any $V$ or $(R-L)$-versal unfolding of an isolated complete intersection singularity is a $T_n$-type mapping (\cite{AVG}).
\end{example}
\subsection{Existence of Lefschetz vanishing cycles}
Let $f:(\CM^m,0) \to (\CM^k,0)$ be a {\em $T_n$-type} mapping. 
Chose a linear mapping $\pi:\CM^k \to \CM^{k-1}$ whose fibre intersect the discriminant at $TA_1^n$ critical
values over some Zariski open subset $U \subset \CM^{k-1}$.
Consider a complex line $L \subset \CM^k$ which is
the preimage of a point in $U$.\\
Take a set of paths connecting a regular value $\e_0 \in (\d S \cap L)$ of $f$
to each of the critical values contained in $L$.\\
Around a point $x$ of type $TA_1^n$, the map $f$
is topologically equivalent to
$$\t f:(x) \mapsto (\sum_{j=1}^{n+1} x_j^2,x_{n+2},\dots,x_k). $$
As $f$ defines a locally trivial fibration above the complement of the discriminant, we may
transport the $n$-dimensional {\em vanishing sphere}
$$S_{x,\e}=\{x \in \RM^k, \t f(x)=(\e,0,\dots,0) \}$$
along the path connecting the critical value $f(x)$ to $\e_0$, we get a sphere in $V$
whose non-oriented isotopy class may a priori depend on the choice of the
critical point.
\begin{lemma}
\label{L::isotopy}
{Two (unoriented) spheres $\a,\b \subset V$ which vanish at two critical points with the same critical value are isotopic.}
\end{lemma}
\begin{proof}
Denote by $ \D \subset M$ the singular locus of the singular fibre
 $ V_{z}=f^{-1}(z) $, $ z \in L \cap \S $ at which the spheres
 $\a,\b$ vanish.\\
The spheres $\a,\b$ vanish at different critical points say
$ x,x' \in V_z $ having the same critical value $z$, i.e.,
$ f(x)=f(x')=z $.\\
In a small neighbourhood $U$ of $x$ the proposition is obviously true.
Cutting a path joining $x$ to $x'$
into small open subsets, so that in each open subset the spheres are isotopic, we get that
the vanishing spheres at $x$ and $x'$ are isotopic after a convenient choice of their orientations.
\end{proof}
\begin{definition}A $T_n$-map germ $f:(\CM^m,0) \to (\CM^k,0)$ is called
{\em calibrated} there exists a differential form $\a \in H^n(V,\CM)$ on the Milnor fibre $V$ of $f$ such that
for any sphere $\g$ vanishing at a $TA_1^n$ critical point, the period integrals $\int_{\g}\a$ is non-zero.
\end{definition}
The calibration allows us to get rid of the ambiguity in the orientation of the vanishing sphere simply by choosing
the orientation for which ${\rm Re\,} \int_{\g}\a>0$ if non-zero and
${\rm \,Im} \int_{\g}\a>0$ otherwise.\\
Therefore to each path connecting a critical value $z$ to the base point $\e_0$ is associated
a well-defined oriented isotopy classes of a sphere, the homology class of
this sphere is called the {\em Lefschetz vanishing cycle} at the critical
point.\\
A set of path connecting the regular value to the critical values of
$f$ contained in $L$ is called {\em weakly distinguished} and the
corresponding vanishing cycles are  called a {\em set of Lefschetz vanishing cycles} associated to these
paths and to the linear mapping $\pi:\CM^k \to \CM^{k-1}$.
The number of cycles contained in such a collection is equal to the
multiplicity of the discriminant at the origin.
\subsection{The free basis theorem}
\begin{definition} A holomorphic map germ $f:(\CM^m,0) \to (\CM^k,0)$
is called {\em simplifiable}\footnote{French word which means that can be simplified.} if at a generic point of the singular locus on the special fibre it has
rank $(k-1)$.
\end{definition}
In the definition {\em generic point} means that the property holds outside an analytic subspace of codimension one
inside the singular locus of the special fibre. 
\begin{theorem}
\label{T::principal}
{Let $f:(\CM^m,0) \to (\CM^k,0)$ be a simplifiable calibrated $T_n$-map germ  with Milnor fibre $V$.
There exists a Zariski open subset $\Omega \subset \Hom_{\CM}(\CM^k,\CM^{k-1})$ such that:
for $\pi \in \Omega$, any set of Lefschetz vanishing cycles
associated to a weakly distinguished basis of paths in $\CM^k$ freely generate the homology group $H_{n}(V)$.}
\end{theorem}
Here and in the sequel, the homology is taken with $\ZM$-coefficient.
\begin{remark} Corollary \ref{C::Lefschetz}
and the Hurewicz isomorphism imply that the groups $H_{n}(V)$ and $\pi_n(V)$ are isomorphic provided that $n>1$.
\end{remark}
\begin{remark}
In the particular case $k=1$, we recover the fact that for isolated hypersurface singularities the Lefschetz
vanishing cycles associated to a morsification of the function generate the group $H_{m-1}(V)$ (\cite{Mil}).
\end{remark}
\begin{corollary}Under the assumptions of the theorem, the multiplicity 
at the origin of the discriminant $m_\S$ of the map-germ $f$ is equal to the
$n$-th Betti number of its Milnor fibre:
$$m_\S={\rm\, rank\ } H_n(V).$$
\end{corollary}
\subsection{Proof of Theorem \ref{T::principal} (first part)}
We prove that the vanishing cycles generate the homology group $H_n(V)$.\\
Denote by $V_0$ the special fibre of a standard representative
$f:X \to S$ above the origin, and by $s$ the dimension of
the singular locus of the variety $V_0$, so that $n=m-k-s$.
According to the results of Subsection \ref{SS::Lefschetz}, we may chose a codimension $s$ linear space
$P=\Ker(u:\CM^m \to \CM^s)$ such that
\begin{enumerate}
\item the pair $(V \cap P, V)$ is
$(m-k-s)$ connected where $s$ denotes the dimension of the singular locus of $V_0$,
\item the variety $V_0 \cap P$
is a reduced complete intersection with an isolated singular point at the origin.
\end{enumerate}
The exact sequence of the pair $(V \cap P,V)$ implies that
the inclusion $V \cap P \subset V$ induces a surjective mapping
$$\xi:H_{n}(V \cap P) \to H_{n}(V).$$
Let us investigate the kernel of this map.\\
The vanishing cycles associated to $(f,u)$ 
generate the middle-dimensional homology of the manifold $V \cap P$
(\cite{Hamm,Mil}).\\
These vanishing cycles may be divided into two classes as follows.
Consider the commutative diagram
$$\xymatrix@!C{ X \ar^-{(f,u)}[r] \ar^-f[rd]& S \times \CM^s  \supset \t \S \ar^-p[d] \\
 & S \supset \S}
$$
The discriminant of the map $(f,u)$ has two components. One of the
components $\t \S$ is the preimage of the discriminant $\S$ of $f$ under the
projection on the first factor $p:S \times \CM^s \to S$.
The other component correspond to regular values of $f$ at which the affine
space $\{ u={\rm constant}\}$ is tangent to the Milnor fibre of $f$.
The corresponding cycles are obviously null-homologous in $V$, i.e., they
lie in the kernel of the map $\xi$.\\
At a critical value of $f$ several distinct vanishing spheres of
$(f,u)$ may vanish. Lemma \ref{L::isotopy} implies that these spheres
may be oriented so to define identical homology classes in $V$.\\
Thus, the image of the map $\xi$ is generated by any set of Lefschetz vanishing cycles associated to $f$.
It remains to prove that they do not satisfy any non-trivial relation. 
\subsection{Simplified mappings}
The Zariski open subset $\Omega$ in the statement of Theorem \ref{T::principal}
is defined in the following proposition.
\begin{proposition}
\label{P::Zariski}
Let $f:(\CM^m,0) \to (\CM^k,0)$ be a simplifiable
complete intersection map germ.
Assume that the special fibre of $f$ has a singular locus of codimension $s>0$. Then there exists a Zariski open subset
$\Omega \subset \Hom_{\CM}(\CM^k,\CM^{k-1})$ such that the special fibre
of the map germ $\pi \circ f$ has a singular locus of codimension at least $s+1$ for $\pi \in \Omega$.
\end{proposition}
\begin{proof}
Consider a representative $f:X \to S$. Take a Whitney stratification
of the special fibre obtained
by refining the stratification by the rank of $f$. Let $X_{sing}$
be the open strata contained in the singular locus of the special fibre. At a point of $X_{sing}$
the rank of $f$ equals $k-1$.\\
Denote by $ Z=\{ (\pi,x) \} \subset \Hom_\CM(\CM^k,\CM^{k-1}) \times X_{sing}$  the analytic subvariety for
which the germ of $\pi \circ f$ has a critical point at $x$.
The projections on the first and second factor give a diagram
$$\xymatrix{ Z \subset  \Hom_\CM(\CM^k,\CM^{k-1}) \times X_{sing} \ar[r]^-{p_1}
                    \ar[d]^-{p_2}& X_{sing}\\
                  \Hom_\CM(\CM^k,\CM^{k-1})} $$
As $f$ is simplifiable, the variety $p_2(Z)$ is of positive codimension. For $\pi \notin p_2(Z)$, the map $\pi \circ f$ is regular
along $X_{sing}$, therefore its
critical space intersect its special fibre inside the boundary of $\bar X_{sing}$ which has
codimension $s+1$. This proves the proposition.
\end{proof}
\subsection{Proof of Theorem \ref{T::principal} (second part).}
We choose a linear mapping $\pi:\CM^k \to \CM^{k-1}$ as in Proposition \ref{P::Zariski}.\\
Denote by $A \subset M$ the Milnor fibre of the involutive pyramidal map $\pi \circ f$
so that $V \subset A$ (if $k=1$ we take $A=M$).\\
Let us consider the exact sequence of the couple $(A,V)$:
$$\cdots \to H_j(V) \to H_j(A) \to H_j(A,V) \to \cdots .$$
According to Corollary \ref{C::Lefschetz}, we have the equalities
$$ H_n(A)=0,\ H_{n+1}(A)=0 .$$ 
This gives an isomorphism
$$ H_n(V) \simeq H_{n+1}(A,V).$$
Denote by $ D_1,\dots, D_l \subset (T \cap L)$ closed disks centred at
 the critical values and containing no other critical value of $ f$.\\
The union of the preimages under $f$ of the disks $D_1,\dots,D_l$ are the union of smooth connected
 manifolds, that we denote by $ \hat D_1,\dots,\hat D_l \subset A$.\\
The excision isomorphism gives the decomposition
$$ H_{n+1}(A,V) \simeq \oplus_{i=1}^l H_{n+1}(\hat D_i,V_i)$$
where $ V_i $ is the fibre of $ f $ at an arbitrary chosen point
of the boundary of $D_i$.\\
The isomorphism $H_{n+1}(V)  \simeq \oplus_{i=1}^l
H_{n+1}(\hat D_i,V_i)$ sends each vanishing sphere $\g_i$ to a ball $C_i$
which bounds a vanishing sphere in $V_i$.\\
The ball $C_i$ generates the group $H_{n+1}(\hat D_i,V_i)$ and it
is non-zero since the map is calibrated. This concludes the proof of the theorem.
\subsection{The Picard-Lefschetz formula for $T_n$ type mappings}
\label{SS::Picard_Lefschetz}
We keep the notations of the preceding subsection. The main observation is given by the following proposition.
\begin{proposition}
\label{P::Slodowy}
Under the assumptions of Theorem \ref{T::principal},
there exists a finite group $\mathcal{G}$ which acts on $H_n(V \cap P)$ preserving the intersection pairing
so that $H_n(V \cap P)/\mathcal{G}$ is isomorphic to $H_n(V)$.
\end{proposition}
\begin{proof}
Take two oriented vanishing spheres $\g_1,\g_2$ in $V\cap P$ which define the same homology class in $V$. The spheres $\g_1,\g_2$ 
vanish at critical points $x_0,x_1$ having the same critical value $z=f(x_0)=f(x_1)$.\\
Denote by $V_{z}$ the fibre of $f$ above $z \in S$
and take a $C^\infty$ path inside the singular locus of $V_z$ connecting $x_0$ to
$x_1$:
$$x:[0,1] \to V_z,\ x(0)=x_0,\ x(1)=x_1.$$
Fix points $p_1,p_2,\dots,p_{d-1} \in P,\ d=\dim P$ and denote by $P_t$ the plane passing through $p_1,\dots,p_{d-1}$, $x(t), x(1-t)$.
Transporting horizontally the homology class in $H_n(V \cap P_t)$ we get a transformation of the cycles which exchange $\g_1$ and $\g_2$,
fixes all vanishing cycles which do not vanish at $z$,
and preserves the intersection pairing. This defines a finite group $\Gt$
acting on  $H_n(V \cap P)$. 
Theorem \ref{T::principal} implies that the quotient of $H_n(V \cap P)$ under
the $\Gt$-action is isomorphic to $H_n(V)$.
\end{proof}
The classical Picard-Lefschetz formula states that for a homology class $a \in H_n(V)$ in the Milnor fibre of an isolated
hypersurface singularity, the monodromy
around a Morse critical value is given by
$$h_*a=a+(a.\D)\D$$
where $\D$ is the cycle vanishing at the critical value (\cite{Lefschetz}). For non-isolated singularities such a formula
does not make sense since the two cycles $a$ and $\D$ do not intersect, if they are in general position.\\

\begin{proposition}{Under the assumptions of Theorem \ref{T::principal},
there is a unique map $D$ such that the diagram
$$\xymatrix{H_n(V \cap P) \ar@{->>}[r]^-p \ar[d]^{P.D.}& H_n(V) \ar^-D[d]\\
H^n(V \cap P) & H^n(V) \ar@{_{(}->}[l]_-i}$$
is commutative. Here $p,i$ are the maps induced by the inclusion $V\cap P \subset V$ and $P.D.$ is the Poincar\'e duality.
Moreover, the monodromy $h_*:H_n(V) \to H_n(V)$ around a critical value is given by
the Picard-Lefschetz type formula
$$h_*=Id+ D\D  $$
where $\D$ is the cycle vanishing at the critical value.}
\end{proposition}
\begin{proof}
According to Theorem \ref{T::principal}, the map $p$ is
surjective and the map $i$ is injective.\\
Take a lifting $s:H_n(V) \to H_n(V \cap P)$ of the map $p$, i.e, $p \circ s=Id$. This lifting
induces a dual mapping $s^*:H^n(V \cap P) \to H^n(V)$.\\
According to Proposition \ref{P::Slodowy}, the map $D= s^* \circ (\, \cdot \,) \circ s$
is independent on $s$ and satisfies the Picard-Lefschetz formula stated above.
\end{proof}
\begin{proposition}
The variation operator $Var:H^n(V) \to H_n(V)$ is well-defined and it is an isomorphism.
\end{proposition}
\begin{proof}
We keep the notations used in the proof of the previous proposition.
As $V \cap P$ is the Milnor fibre of an isolated
complete intersection singularity, the variation operator
$\widetilde{Var}: H^n(V \cap P) \to H_n(V \cap P)$ exists (see e.g. \cite{AVGII}).\\
We define the variation operator in $H^n(V)$ using the commutative diagram
$$\xymatrix{H^n(V) \ar@{^{(}->}[r]^-i \ar[d]^{Var}& H^n(V \cap P) \ar^{\widetilde{Var}}[d] \\
H_n(V) & H_n(V \cap P) \ar@{->>}[l]^p }$$
According to Theorem \ref{T::principal}, we may chose a basis of vanishing
cycles which generate the group $ H_{n}(V)$ freely.
Take for  $H^{n}(V)$ the dual basis to this basis. In these basis, the Picard-Lefschetz formula implies that the variation operator is triangular with
  diagonal elements equal to $1$ or $-1$ (see e.g. \cite{AVGII},
  Chapitre I, section 2.5). This proves the proposition.
\end{proof}

\section{Involutive mappings in symplectic spaces}
\subsection{Pyramidal mappings}
Let us now consider the case, where $\CM^m=\CM^{2n}$ is endowed with the holomorphic symplectic form
$\omega=\sum_{i=1}^n dq_i \w dp_i$. The symplectic form induces a Poisson Bracket
$$\{ f , g \}\omega^n=df \w dg \w \omega^{n-1} $$
which leads to the standard formula
$$\{ f , g \}=\sum_{i=1}^n \d_{q_i}f \d_{p_i}g- \d_{p_i}f \d_{q_i}g .$$
The interior product with $\omega$ induces an isomorphism between the sheaves of
holomorphic one-forms and that of holomorphic vector fields. Given a holomorphic function
$H:U \to \CM$, defined on an open subset
$U \subset \CM^{2n}$, the vector field associated to the holomorphic one-form $dH$ is called the {\em Hamiltonian
vector field} associated to $H$.
\begin{definition}
A holomorphic map $f:X \to S,\ S \subset \CM^k$ is called
{\em involutive} if the Poisson brackets of its components vanish.
\end{definition}
\begin{definition}A {\em pyramidal map}
$f=(f_1,\dots,f_k):U \to \CM^k$ where $U$ is an open subset in $\CM^{2n}$ 
is a stratified involutive map such that for any critical point $x$ lying on a strata $X_\a$,
the Hamiltonian vector fields $v_1,\dots,v_k$ give local coordinates around $x$ on $f^{-1}(x) \cap X_{\a}$, i.e.:
$$T_x(f^{-1}(x) \cap X_\a)={\rm Span}\{ v_1(x),\dots,v_k(x) \}$$
at any point $x$ of a strata $S$ inside the singular locus of a fibre.
\end{definition}
A similar notion holds for germs, in that case it is sufficient to check
the property on the special fibre.
\begin{example}For $k=1$, a map is pyramidal provided it has isolated critical points.
\end{example}
\begin{example}
Let $f=(f_1,f_2):(\CM^4,0) \to (\CM^2,0)$ be an involutive complete
intersection map germ. Assume that some linear combination of $f_1$ and $f_2$ has an isolated critical point at the origin, then $f$ is pyramidal
\end{example}

\subsection{Pyramidal involutive mappings satisfy the $a_f$-condition}
\begin{proposition}
\label{P::standard} Any involutive pyramidal map germ
$f:(\CM^{2n},0) \to (\CM^k,0)$ on a symplectic space satisfies the strong $a_f$-condition.
In particular, the restriction of the map $f$ above the complement
of the discriminant defines a $C^\infty$ locally trivial fibration.
\end{proposition}
\begin{proof}
We use the notations of Definition 1.\\
Let $f:X \to S$ be a representative of the germ $f$ satisfying
the first two conditions of Definition 1 and denote by $X_s$ the fibre above $s \in S$.\\
Chose a Whitney stratification of $X$ obtained by  a refinement of the stratification by the rank of $f$
and whose intersection with the special fibre defines a Whitney stratification.\\
Let $x$ be a singular point of the fibre $X_0$, the pyramidality assumption
implies that the tangent space to the strata of $x$ is
generated by the Hamiltonian vector fields associated to $f_1,\dots,f_k$.
These Hamiltonian fields are defined in a neighbourhood
of $x$ and therefore by transversality the tangent space $T_x S$ at $x$ is the limit of the tangent spaces 
to the adjacent fibres.\\
(Take a vector $v_x$ tangent at $x$ to
the strata of $x$. The pyramidality assumption implies that there
exist $a_1,\dots,a_n \in \CM$ such that $v_x=\sum a_i v_i(x)$ where
$v_i$ is the Hamiltonian field associated to $f_i$.
The vector field $ \sum a_i v_i$ is globally defined on $M$ and it is
tangent to the fibres of $f$.)
This proves the proposition.
\end{proof}

\subsection{Pyramidal mappings in symplectic spaces are calibrated}
The proof of the following proposition is straightforward.
\begin{proposition}
\label{P::produit}
Let $f:(\CM^{2n},0) \to (\CM^k,0)$ be the germ of an involutive map with a
non-vanishing differential at the origin.
There exist biholomorphic maps $\p:(\CM^{2n},0) \to
(\CM^{2n},0)$, $\psi:(\CM^k,0) \to (\CM^k,0)$ such that $ \p$ is
symplectic and
$$\psi \circ f \circ \p=(p_1,\dots,p_j,g)$$
where $j$ denotes the rank of $df(0)$.
\end{proposition}
The symplectic form $\omega$ being closed, we may chose a differential one-form
$\a$ such that $d\a=\omega$. In local coordinates, we have $\a=\sum_{i=1}^n p_idq_i$.
\begin{proposition}Any pyramidal involutive $T_n$-mapping is calibrated with respect to the differential
form $\a \w (d\a)^{n-k}$.
\end{proposition}
\begin{proof}
According to the previous proposition at any smooth point of a fibre, we may chose local Darboux
coordinates such that $f(q,p)=(p_1,\dots,p_k)$. These choice of coordinates shows that the restriction of the differential form
$\a \w (d\a)^{n-k}$ to the fibres of $f$ is closed and hence define a
cohomology class.\\
Any vanishing sphere $\g$ in the Milnor fibre is the boundary of a ball $B$. By Stokes formula, we get
$$\int_{\g}\a \w (d\a)^{n-k}=\int_{B} (d\a)^{n-k+1} \neq 0$$
This proves the proposition.
\end{proof}
\subsection{Arnold-Liouville manifolds}
Consider the pyramidal involutive map 
$$f:(\CM^6,0) \to
    (\CM^2,0),\ (q_1,q_2,q_3,p_1,p_2,p_3) \mapsto
(p_1q_1+p_2q_2, \ p_2q_2+p_3q_3).$$ The discriminant of $f$
 consists in the germ at the origin of  the three lines
$\{ s_1=0 \} \cup \{s_2=0 \} \cup \{
s_1=s_2 \}$. Therefore, the first three homotopy groups of the Milnor
 fibre $V$ of $f$ are $\pi_1(V)=\pi_2(V)=0$ and $\pi_3(V) =H_3(V)= \ZM^3$. So, the
complex manifold $V$ is a 2-connected affine
$4$-fold with non-vanishing three dimensional homology group.\\
This example can be generalised as follows.
Let $g=(g_1,\dots,g_n):(\CM^{2n},0) \to (\CM^n,0)$ be the map germ
defined by  $g_i=p_i q_i$. Denote by $R:\CM^n \to \CM^k$ a surjective
linear map such that $R \circ g$ is pyramidal.
The Milnor fibre of the pyramidal involutive map  $f=R \circ g:(\CM^{2n},0)
\to (\CM^k,0)$ is a $2(n-k)$-connected $2n-k$-dimensional
complex affine manifold whose Betti number
$b_{2(n-k)+1}$ equals the binomial coefficient
$\begin{pmatrix} n  \\ k-1 \end{pmatrix}$. In case $k=n$, the fibres are homotopic to Arnold-Liouville
tori, therefore we call such manifolds arising from integrable system {\em Arnold-Liouville manifolds}.

\subsection{The Henon-Heiles integrable system}
\begin{proposition} The Milnor fibre $V$
of the Henon-Heiles integrable system  $H=(H_1,H_2):(\CM^4,0) \to (\CM^2,0)$
with
$$H_1=p_1^2+p_2^2-4q_2^3-2q_1^2q_2,\ H_2=q_1^4-4q_1^2q_2^2+4p_1(q_1p_2-q_2p_1)  $$
has its first Betti number equal to 4.
\end{proposition}
\begin{proof}
The discriminant of the Henon-Heiles integrable system
$$H_1=p_1^2+p_2^2-4q_2^3-2q_1^2q_2,\ H_2=q_1^4-4q_1^2q_2^2+4p_1(q_1p_2-q_2p_1)  $$
has two components. One of the component is the curve $s_2^3=s_1^4$ the other component is the straight line
$\{ s_2=0 \}$(\cite{Gerardy}).\\
The fibres of the map $f=(H_1,H_2):\CM^4 \to \CM^2$ are open parts of Prym varieties associated to the two
fold covering of an elliptic curve by a hyperelliptic genus 3 curve (\cite{Gerardy,Lesfari}).
Therefore the singular locus of a generic singular fibre is connected, this can also be seen directly: a direct computation shows that
the critical locus of $f$ has two components, each component intersecting the special fibre along a reduced connected curve.
Applying Theorem \ref{T::principal}, we get that the first
Betti number of the Milnor fibre is equal to 4. This proves the proposition.
\end{proof}
\section{The vanishing topology of coadjoint orbits}
\subsection{The Steinberg mapping of a semi-simple Lie algebra}
A {\em Poisson structure} in $\CM^n$ is an antisymmetric biderivation
$$\{ \cdot , \cdot \}:\OM_{\CM^n} \times \OM_{\CM^n} \to \OM_{\CM^n}$$
satisfying the Jacobi identity.\\
A Poisson structure on a finite dimensional vector space $V$ is given by specifying the values of the Poisson
bracket on the space of linear forms $V^*$. By identifying a Lie algebra with its bidual we see that the Lie Bracket defines a Poisson bracket on the dual of the Lie algebra
$$\{ a,b \}:=[a,b]$$
called the {\em Kirillov-Kostant-Souriau Poisson structure}. This Poisson structure defines a symplectic structure on coadjoint orbits, that is, identifying $\gk^*$ with the image of a faithful representation,
the orbits under conjugacy.\\
The {\em Steinberg mapping} $S: \gk \to \hk/W$ is the map sending an element $x=x_s+x_n$
to the class of its semi-simple part (in the Jordan decomposition) modulo the action of the Weyl group.
By Chevalley's theorem there is a scheme isomorphism of the Cartan subalgebra
modulo the Weyl group and the adjoint quotient
$$\hk/W \approx \gk/G.$$
By Cartan's theorem, as $\mathfrak{g}$ is simple,
the Killing form $\gk$ is non-degenerate and therefore coadjoint orbits and 
adjoint orbits can be identified. 
\begin{proposition}
\label{P::Steinberg}
The germ of the Steinberg mapping at any point is a simplifiable
calibrated $T_2$-mapping.
\end{proposition}
\subsection{Some remarks concerning the vanishing topology of coadjoint orbits}
Proposition \ref{P::Steinberg}
enables us to apply Theorem \ref{T::principal} for the Steinberg mapping.
Remark that in particular the Steinberg mapping defines a $C^\infty$ locally trivial fibration above the complement of its discriminant.\\
We combine this result with the following theorem due to Brieskorn.
\begin{theorem}[\cite{Brieskorn_71}]
Denote by $(m_{ij})_{i,j=1,\dots,r}$
the Coxeter matrix associated to a semi-simple Lie algebra
$\mathfrak g$. Then, there exists homotopy classes $b_1,\dots,b_r$
associated to a basis of Lefschetz vanishing cycles such that
the group $\pi_1(\hk/W \setminus \S) $ admits a presentation with
generators $b_1,\dots,b_r$ and relations
$$\underbrace{b_ib_jb_i \dots}_{m_{ij}\, times}= \underbrace{b_jb_ib_j \dots}_{m_{ij}\, times}\ {\rm whenever}\ m_{ij} \geq 3$$
\end{theorem}
By adding the relations $b_i^2=1$ to the group $\pi_1(\hk/W \setminus \S) $, we get the {\em Coxeter group} $G$
associated to the Lie algebra.
The group $G$
admits a faithful representation as a finite reflexion group in $\RM^n$. The induced
representation on $\pi_1(\hk/W \setminus \S) $ is called {\em the Coxeter representation}.\\
According to the standard Picard Lefschetz formula and to its version for non-isolated singularities
(see Subsection \ref{SS::Picard_Lefschetz}), the monodromy of Lefschetz vanishing cycles
is a reflection in the hyperplane orthogonal to the vanishing cycle. By Theorem \ref{T::principal}, all these hyperplane are distinct
and Brieskorn's theorem gives the angle between this hyperplanes, therefore we get
the following result.
\begin{theorem}[Compare \cite{Rossmann}] The monodromy representation
$$\pi_1(\hk/W \setminus \S) \to GL(r,H_2(V) \otimes\RM),$$
 associated to the germ of the Steinberg mapping of a simple Lie algebra at a critical point $x$ with $S(x)=0$,
coincides with the Coxeter representation.
\end{theorem}
Applying the results of Subsection \ref{SS::Picard_Lefschetz} and
comparing this Theorem with the results of Brieskorn, Esnault and Slodowy (\cite{Brieskorn_ICM,Esnault,Slodowy}),
we get that, at a subregular point,
for the simple Lie algebras of type $A,D,E$ there is an isomorphism
$H_2(V \cap P) \approx H_2(V)$ and
that the group $\mathcal G$ is trivial, while for the ones of type $B,C,F,G$
the group $\mathcal{G}$ are respectively the group $\ZM/2\ZM$ for $B,C,F$ and
the permutation group ${\mathcal S_3}$ for the group $G_2$. These groups attached to the $B,C,F,G$ Lie algebras were discovered independently
by Arnold and Slodowy
(\cite{Arnold_BCF} and \cite{Slodowy},Chapter 6.2).
\subsection{The $a_f$-condition and Casimir functions}
Recall that a holomorphic function is called a {\em Casimir function} if it Poisson commutes with any other holomorphic function.
\begin{proposition}
Let $f:(\CM^n,0) \to (\CM^k,0)$ be a complete intersection
map germ such that  the components of $f$ are Casimir functions,
then $f$ satisfies the strong $a_f$-condition.
\end{proposition}
\begin{proof}
The proof is similar to that of pyramidal involutive mappings as any vector $v_x$ tangent to a strata
of the special fibre at a point $x$ extends to a Hamiltonian vector field in a neighbourhood of $x$.
Indeed, choose a function $h:U \to \CM$ such that $v_x$ is the Hamiltonian vector field $v$
of $h$ evaluated at $x$ in some neighbourhood $U$ of the point $x$.
As the components of $f$ Poisson commutes with $h$, the vector field $v$ is tangent to the fibres of $f$.
\end{proof}
Applied to the Steinberg mapping the proposition shows that it defines a locally trivial fibration
in any small open ball containing the origin.
\subsection{Proof of Proposition \ref{P::Steinberg}}
An element $x \in \gk$ is called {\em regular} if its
centraliser has dimension $r+2$ where $r$ is the rank of the Lie algebra, i.e,
the dimension of a Cartan subalgebra. A theorem due to Kostant and Steinberg states that the non-regular elements are the critical
values of the Steinberg mapping (\cite{Kostant,Steinberg}).\\
An element $x \in \gk$ is called {\em subregular} if its
centraliser has dimension $r+2$.
Any singular fibre of the Steinberg
mapping contains a subset of subregular elements which is dense in the singular locus of the fibre (\cite{Steinberg}).
The rank of the Steinberg mapping at such point equals $r-1$ (\cite{Brieskorn_ICM,Slodowy}).
\begin{lemma} Let $x \in \gk$ be a critical point of the Steinberg mapping.
If $S(x)$ is a smooth point of the discriminant then is $x$ semisimple and subregular.
\end{lemma}
\begin{proof}
We consider the Jordan decomposition $x=x_s+x_n$ it into semisimple and nilpotent part and
 use the criterion (\cite{Steinberg}):\\
$x$ is subregular in $ \gk \iff x_n$ is subregular in $C(x_s)$.\\
Here $C(x_s)$ denotes the centraliser of $x_s$.\\
Chose a Cartan subalgebra $ \hk$ containing $x_s$ and let $\Phi$ be the associated root system.
The root space decomposition
$$ \gk=\hk\oplus \bigoplus_{\a \in \Phi}  \gk_\a$$
shows that $C(x_s)=\hk\oplus \gk_\b \oplus \gk_{-\b}$
where $\a$ is the root which annihilates $x_s$; thus $x_s$ is subregular. As $x_n$ is subregular in $C(x_s)$
and $C(x_s)$ has dimension $r+2$ we get that $C(x_n)=C(x_s)$; in particular
$x_n$ commutes to all element in the Cartan subalgebra $\hk$, thus $x_n \in \hk$ and consequently $x_n=0$.
\end{proof}
The lemma implies in particular that the singular locus of $S^{-1}(S(x))$ is connected for $x$ in the smooth locus of the discriminant.
\begin{lemma}{\rm (Compare \cite{Arnold_matrices})}{\ If $x$ is a subregular semisimple element of the rank $r$ Lie algebra $\gk$,
then there exists a commutative diagram
of holomorphic map germs
$$\xymatrix{(\gk,x) \ar[r]^-S \ar[d]^-\p & (\hk/W,S(x)) \ar[d]^-\psi \\
             (\sl(2,\CM),0) \times (\CM^{r-1},0) \ar[r]^-{(det,Id)} & (\CM,0) \times (\CM^{r-1},0)}$$}
\end{lemma}
\begin{proof}
Like in the previous lemma, we consider the Chevalley-Cartan decomposition
$$ \gk=\hk \oplus \bigoplus_{\a \in \Phi'} \gk_\a \oplus \gk_\b
\oplus \gk_{-\b},\ \Phi'=\Phi \setminus \{ \b,-\b \}.$$
so that the roots $\pm \b$ annihilate $\ x \in \hk$.\\
The tangent space to the adjoint orbit of $x$ is
$$[x,\gk]= \bigoplus_{\a \in \Phi'} \gk_\a$$
therefore the vector space $\hk \oplus \gk_\b \oplus \gk_{-\b}$
is transversal to it.
Consequently, the germ at $x$ of the Steinberg mapping factorises
$$\xymatrix{(\gk,x) \ar[r]^-S \ar[d] & (\hk/W,S(x)) \\
             (\hk\oplus \gk_\b \oplus \gk_{-\b},x) \ar[ru]^-f}.$$
Consider the hyperplane $\hk_1=Ker \b$ and consider a supplementary
one dimension vector space to it $\hk_2 \subset \hk$ , the rank one
Lie algebra $\hk_2\oplus \gk_\b \oplus \gk_{-\b}$
is isomorphic to ${\sl(2,\CM)}$, therefore we get an isomorphism
of Lie algebras
$$\hk\oplus \gk_\b \oplus \gk_{-\b} \approx
\sl(2,\CM) \oplus \hk_1.$$
Via this identification, the Lie subgroup $H \subset G$ which fix the
transversal act on $\sl(2,\CM)$ by conjugation. Moreover, by Chevalley's
theorem the action of $G$ on $\hk_1$ is induced by that of the Weyl group
and by Kostant-Steinberg characterisation of regular elements, the Steinberg mapping
restricted to $\hk_1$ has maximal rank at $x$.
This concludes the proof of the lemma.
\end{proof} 
{\em Acknowledgements.}{ The author thanks D.T. L{\^e}, E. Lettelier, J. Pagnon,
F. Pham, B. Teissier, D. van Straten for fruitful discussions
and V.I. Arnold for his teaching.}

\bibliographystyle{amsplain}
\bibliography{master}
\end{document}